\documentclass{amsart}
\usepackage{amssymb,latexsym,amsmath,amsthm,enumitem,mathrsfs,geometry}
\usepackage{fullpage}
\usepackage{fancyhdr}
\usepackage{color}
\usepackage{stmaryrd}
\usepackage{mathrsfs}
\usepackage{hyperref}
\usepackage{lineno}
\usepackage{diagbox}
\usepackage{array}
\usepackage{tikz}
\usetikzlibrary{arrows}
\usetikzlibrary{positioning}
\usepackage{float}
\usepackage{bbm}
\usepackage{cleveref}
\usepackage{comment}

\allowdisplaybreaks

\newcommand{\1}{\mathbbm{1}}
\newtheorem{theorem}{Theorem}
\newtheorem{lem}{Lemma}

\newtheorem*{theorem*}{Theorem}
\definecolor{pink}{rgb}{1,.2,.6}
\definecolor{orange}{rgb}{0.7,0.3,0}
\definecolor{blue}{rgb}{.2,.6,.75}
\definecolor{green}{rgb}{.4,.7,.4}
\definecolor{purple}{RGB}{127,0,255}

\makeatletter
\@namedef{subjclassname@2020}{\textup{2020} Mathematics Subject Classification}
\makeatother

\begin{document}

\title{On a smoothed average of the number of Goldbach representations}

\author[Goldston]{Daniel A. Goldston}
\address{Department of Mathematics and Statistics, San Jose State University}
\email{daniel.goldston@sjsu.edu}

\author[Suriajaya]{Ade Irma Suriajaya$^{*}$}
\address{Faculty of Mathematics, Kyushu University}
\email{adeirmasuriajaya@math.kyushu-u.ac.jp}
\thanks{$^{*}$ The second author was supported by JSPS KAKENHI Grant Numbers 18K13400 and 22K13895, and also by MEXT Initiative for Realizing Diversity in the Research Environment.}
\date{\today}

\dedicatory{Dedicated to the memory of Eduard Wirsing}
\begin{abstract}
Assuming the Generalized Riemann Hypothesis for the zeros of the Dirichlet $L$-functions with characters modulo $q$, we obtain a smoothed version of the average number of Goldbach representations for numbers which are multiples of a positive integer $q$. Such an average was first considered by Granville \cite{Gran07,Gran08} but without any smoothing factor.
In this short article, we also show how the smoothing can be removed. 

\end{abstract}

\maketitle


\section{Introduction and Statement of Results}

In this paper all sums are over positive integers unless otherwise indicated. Let
\begin{equation*}
\psi_2(n) = \sum_{m+m'=n} \Lambda(m)\Lambda(m'),
\end{equation*}
where $\Lambda$ is the von Mangoldt function, defined by $\Lambda(n)=\log p$ if $n=p^m$, $p$ a prime and $m\ge 1$, and $\Lambda(n)= 0$ otherwise. Thus $\psi_2(n)$ counts the \lq \lq Goldbach" representations of $n$ as sums of both primes and prime powers, and these primes are weighted to make them have a \lq\lq density" of 1 on the integers. 

The average number of Goldbach representations was first studied by Landau \cite{Landau1900} in 1900. Letting $\pi_2(n) =\sum_{p+p'=n}1$ with $p$ and $p'$ primes, he proved using the prime number theorem with error term that
\begin{equation*} \sum_{n\le N} \pi_2(n) = \frac{N^2}{2 \log^2\!N} + O\left( \frac{N^2}{\log^3\!N}\right). \end{equation*}
Landau's result implies
\begin{equation} \label{Landau}
G(N) \sim \frac{ N^2}{2},\qquad \text{where} \qquad G(N) := \sum_{n\le N}\psi_2(n).
\end{equation}
Hardy and Littlewood \cite{HardyLittlewood1919,HardyLittlewood1922} used this result in testing conjectures concerning $\psi_2(n)$. More recently, 
Fujii \cite{Fujii1,Fujii2,Fujii3} in 1991 proved the following theorem concerning the average number of Goldbach representations, which for the first time directly connected the error term in Landau's formula to the zeros of the Riemann zeta-function and the Riemann Hypothesis (RH). 
\begin{theorem*}[Fujii] Assuming the Riemann Hypothesis, so that the complex zeros $\rho=\beta +i\gamma$ of the Riemann zeta-function $\zeta(s)$ have $\beta = 1/2$. Then
\begin{equation} \label{FujiiThm} G(N) = \frac{N^2}{2} - 2\sum_{\rho} \frac{N^{\rho+1} }{\rho(\rho+1)} + O(N^{\frac43}(\log N)^{\frac43}). \end{equation}
\end{theorem*}

The error term above was improved by Bhowmik and Schlage-Puchta in \cite{BhowPuchta2010} to $O(N\log^5\!N)$ assuming RH, and this was refined by Languasco and Zaccagnini \cite{Lang-Zac1} who obtained the following result. 
\begin{theorem*}[Languasco-Zaccagnini] Assuming the Riemann Hypothesis, we have
\begin{equation} \label{LZ} G(N) = \frac{N^2}{2} - 2\sum_{\rho} \frac{N^{\rho+1}}{\rho(\rho +1)} + O(N\log^3\!N). \end{equation}
\end{theorem*}
\noindent
As shown in \cite{Gold-Yang}, the above asymptotic can also be obtained by extending the method of \cite{BhowPuchta2010}.
Bhowmik and Schlage-Puchta in the previously mentioned paper \cite{BhowPuchta2010} also proved unconditionally that
\begin{equation*} 
G(N) = \frac{N^2}{2} - 2\sum_{\rho} \frac{N^{\rho+1}}{\rho(\rho+1)} + \Omega(N\log\log N),
\end{equation*}
and therefore the error term in \eqref{LZ} is close to best possible. 

The Riemann Hypothesis was used in Fujii's Theorem to handle the error term, and this obscures the dependence of the error term on the zeros. This has been addressed in the recent papers \cite{BhowRuzsa2018} and \cite{{Bhowmik-H-M-SGoldbach2019}} which obtain an equivalence between the size of the error in \eqref{Landau} with a zero-free region for zeros of the Riemann zeta-function. 

Granville \cite{Gran07,Gran08} introduced the average number of Goldbach representations for integers which are multiples of $q$ by defining\footnote[2]{We have $G_q(N)=0$ if $q>N$. }
\begin{equation} \label{Gq} G_q(N) := \sum_{\substack{n\le N\\ q\mid n}}\psi_2(n), \qquad \text{for} \qquad 2\le q\le N. \end{equation}
When $q=1$ we write $G(N)=G_1(N)$ for the case when we average over all Goldbach representations. Granville proposed a remarkable asymptotic formula for $G_q(N)$, but at the time the error term in \eqref{FujiiThm} appeared to interfere with the conjectured error of the proposed formula. This problem was removed by Bhowmik and Schlage-Puchta \cite{BhowPuchta2010} and now the result \eqref{LZ} is available. Thus we have, assuming the Generalized Riemann Hypothesis (GRH) holds for all Dirichlet $L$-functions $L(s,\chi)$, with characters $\chi({\rm mod}\, q)$, that 
\begin{equation} \label{Gqformula} G_q(N) = \frac{1}{\phi(q)}G(N) + O(N\log^C\!N) , \end{equation}
for some fixed constant $C$. 
This can be obtained by the method of \cite{Lang-Zac1} and \cite{Suzuki2017}, or the methods of \cite{Gran07} and \cite{Bhowmik-H-M-SGoldbach2019}. 
The more general problem of considering averages of $\psi_2(n)$ in arithmetic progressions, or having each prime in the Goldbach representation of $n$ in its own arithmetic progression has also been examined in \cite{Suzuki2017} and \cite{{Bhowmik-H-M-SGoldbach2019}}.

Our goal in this note is to study a smoothed average for the number of Goldbach representations which are multiples of $q$. 
We define this average to be
\begin{equation} \label{FqN}
F_q(N) = \sum_{\substack{n\\ q\mid n}}\psi_2(n)e^{-n/N},
\end{equation}
and when $q=1$ we let $F(N) := F_1(N)$ for averaging over all positive integers. By smoothing with this power series weight, we are able to obtain for $F_q(N)$ and $F(N)$ results corresponding to \eqref{Gqformula} and \eqref{LZ} for $G_q(N)$ and $G(N)$ with minimal effort.

\begin{theorem} \label{thm-smoothed}
Suppose the Generalized Riemann Hypothesis holds for the Dirichlet $L$-functions $L(s,\chi)$ with characters $\chi({\rm mod}\, q)$. Then for $q\ge 2$
\begin{equation} \label{thm2a} F_q(N) = \frac{1}{\phi(q)}F(N) + O\left(N\log N \log q \right) , \end{equation}
and when $q=1$
\begin{equation} \label{thm2b} F(N) = N^2 - 2\sum_{\rho} \Gamma(\rho)N^{\rho+1} + O(N),
\end{equation}
where the sum runs over the non-trivial zeros of the Riemann zeta-function, and $\Gamma(s)$ is the gamma function.
\end{theorem}
\noindent
The second equation actually depends on RH, which is included from our assumption on $L(s,\chi_0)$.

The method we use is a combination of the methods of the papers \cite{GS21} and \cite{FGIS21}. We can use this method to prove \eqref{Gqformula} but the proof is much more complicated than the proof of Theorem \ref{thm-smoothed}. 
This paper is organized as follows. In Section \ref{sec2} we set up the results needed in the proof of Theorem \ref{thm-smoothed}. These are actually proved in greater generality than necessary so that they also apply for $G_q(N)$. In Section \ref{sec3} we prove Theorem \ref{thm-smoothed}. Finally, in Section \ref{sec4} we outline how to apply this method to $G_q(N)$.

\section{Lemmas}
\label{sec2}

Let 
\begin{equation} \label{F_q} F_q(z) := \sum_{\substack{n\\ q\mid n}}\psi_2(n)z^n,\end{equation}
where we take, for $N\ge 4$,
\begin{equation} \label{zrNalpha} z=r e(\alpha), \quad r = e^{-1/N} , \quad \text{and} \quad e(\alpha) = e^{2\pi i\alpha}.
\end{equation}
We also have, on letting $F(z)= F_1(z)$, that
\begin{equation}\label{F=Psi^2} F(z) = \Psi(z)^2, \qquad \text{where} \qquad \Psi(z) := \sum_{n }\Lambda(n)z^n.\end{equation}

To evaluate $F_q(z)$, we use
\begin{equation*} \Psi(z,\chi) := \sum_n \chi(n)\Lambda(n) z^n, \end{equation*}
where $\chi$ is a Dirichlet character $({\rm mod}\, q)$. 
\begin{lem} \label{Fqlem}
For $N\ge 4$ and $q\ge 2$, we have 
\begin{equation} \label{lemma1result} F_q(z) = \frac{1}{\phi(q)}\sum_{\chi ({\rm mod}\, q)} \chi(-1)\Psi(z,\chi)\Psi(z,\overline{\chi}) + O\left((\log N\log q)^2\right).
\end{equation}
If $q=1$ we have the identity \eqref{F=Psi^2}.
\end{lem}

\noindent\textit{\textbf{Proof} (Proof of Lemma \ref{Fqlem})}
For $N\ge 4$ and $q\ge 2$, we have 
\begin{equation}\begin{split}\label{S1S2}
F_q(z) &= \sum_{q|m+m'} \Lambda(m)\Lambda(m')z^{m+m'}
= \sum_{m+m' \equiv 0 \,({\rm mod}\, q)} \Lambda(m)\Lambda(m')z^{m+m'} \\
&= \left(\sum_{\substack{m\\ (m,q)=1}}+\sum_{\substack{m\\ (m,q)>1}} \right)\left(\sum_{\substack{m' \\ m' \equiv -m \,({\rm mod}\, q)}} \Lambda(m)\Lambda(m')z^{m+m'}\right)
 \\
&=: S_1 +S_2.\end{split}
\end{equation}

We introduce characters into $S_1$ using the relation, for $(a,q)=1$,
\begin{equation*} \frac{1}{\phi(q)} \sum_{\chi ({\rm mod}\, q)}\chi(a)\overline{\chi}(n) = \1_{n\equiv a ({\rm mod}\, q)}.
\end{equation*}
Thus
\begin{align*}
S_1
&= \sum_{m,m'} \left( \frac1{\varphi(q)} \sum_{\chi ({\rm mod}\, q)} \chi(-m)\overline{\chi}(m') \right) \Lambda(m)\Lambda(m')z^{m+m'} \\
&=\frac1{\varphi(q)} \sum_{\chi ({\rm mod}\, q)} \chi(-1) \left( \sum_m \chi(m)\Lambda(m)z^m \right) \left(\sum_{m'} \overline{\chi}(m')\Lambda(m')z^{m'} \right) \\
&= \frac1{\varphi(q)} \sum_{\chi ({\rm mod}\, q)} \chi(-1) \Psi(z, \chi)\Psi(z, \overline{\chi}).
\end{align*}

For $S_2$, we note that $(m',q)>1$, and therefore
\[ S_2 \ll \left(\sum_{\substack{m\\ (m,q)>1}}\Lambda(m)r^m\right)^2.\]
We have 
\begin{equation} \label{Psiqtrivial} \begin{split}\sum_{\substack{m\\ (m,q)>1}}\Lambda(m)r^m &= \sum_k \sum_{p|q} (\log p)r^{p^k} \\ &
\le \left( \sum_{p|q} \log p \right)\left( \sum_k e^{- \frac{2^k}{N}}\right) \\&
\le \log q \left( \sum_{2^k\le N} 1 + \sum_{2^k> N} e^{- \frac{2^k}{N}}\right) \\&
\le \log q \left( \left \lfloor \frac{\log N}{\log 2}\right \rfloor + \sum_{j} e^{- j}\right)\\& \ll \log N \log q, \end{split} \end{equation}
and \eqref{lemma1result} follows.

If $q=1$, we see in \eqref{S1S2} that $S_2=0$ and
\[F(z) = \sum_{m,m'}\Lambda(m)\Lambda(m')z^{m+m'}= \Psi(z)^2.\]

\qed

Let $\chi_0$ denote the principal character, and define 
\begin{equation*} 
E_0(\chi) :=
\begin{cases}
1, & \text{if $\chi=\chi_0$,}\\
0, & \text{otherwise.}
\end{cases}
\end{equation*}

\begin{lem} \label{Psi(z)Explicit}
Let $z= e^{-w}$, so that from \eqref{zrNalpha} we have $w = \frac{1}{N} - 2\pi i\alpha$. Then for $q\ge 1$ 
\begin{equation} \label{PsiExplicit} \Psi(z,\chi) = \sum_n \chi(n)\Lambda(n) e^{-nw} = \frac{E_0(\chi)}{w} -\sum_{\rho}\Gamma(\rho)w^{-\rho} + O(\log(2qN)\log{2q}), \end{equation}
where the sum is over the non-trivial zeros of $L(s,\chi)$. In particular, if $w=\frac{1}{N}$ then for $\chi\neq \chi_0$ and $q\ge 3$ we have
\begin{equation} \label{Psi(r)explicit} \Psi(r,\chi) = \sum_n \chi(n)\Lambda(n) e^{-n/N} = -\sum_{\rho}\Gamma(\rho)N^{\rho} + O(\log(qN)\log{q}), \end{equation}
and for $\chi=\chi_0$ and $q\ge 2$ 
\begin{equation} \label{PsiChi_0} \Psi(z,\chi_0) = \sum_{\substack{ n\\ (n,q)=1}} \Lambda(n)z^n =\Psi(z) + O(\log N \log q).\end{equation} 
Finally, we have
\begin{equation} \label{PsiNexplicit} \Psi(r) = \sum_n \Lambda(n) e^{-n/N} = N - \sum_{\rho}\Gamma(\rho)N^{\rho} -\log 2\pi+ O\left(\frac1{N}\right), \end{equation}
with the sum over the non-trivial zeros of the Riemann zeta-function.
\end{lem}
The proof of \eqref{PsiExplicit} may be found in \cite[Lemma 2]{Lang-Zac2}. It is stated there that GRH is assumed, but the proof holds unconditionally. See also \cite[Lemma 2.1]{Suzuki2017}. In this paper we only need the special cases where $z=r$. We see \eqref{PsiChi_0} follows immediately from \eqref{Psiqtrivial}. 
In the explicit formula for $\Psi(r)$ \eqref{PsiNexplicit}, the error term is from \cite[12.1.1 Exercise 8(c)]{MontgomeryVaughan2007} which follows from the original proof in \cite{HardyLittlewood1918}.

\section{Proof of Theorem \ref{thm-smoothed}}
\label{sec3}
We take $z=r=e^{-1/N}$ in \eqref{F_q} and have
\[ F_q(r) = \sum_{\substack{n\\ q\mid n}}\psi_2(n)e^{-n/N} = F_q(N), \]
using the notation \eqref{FqN}.
By Lemma \ref{Fqlem} we have for $q=1$
\begin{equation} \label{Fq=1} F(N) = \Psi(r)^2,\end{equation}
and for $q\ge 2$,
\begin{equation} \label{Fqformula} F_q(N) = \frac{1}{\phi(q)}\sum_{\chi ({\rm mod}\, q)} \chi(-1)|\Psi(r,\chi)|^2 + O\left((\log N\log q)^2\right). \end{equation}
In Lemma \ref{Psi(z)Explicit} we use $\Gamma(\rho) \ll e^{-|\gamma|}$ and the estimate for all $t\in\mathbb{R}$ 
\[ \sum_{t<\gamma \le t+1} 1\ll \log (2q(|t|+1)) \]
 to see the sum over zeros in \eqref{Psi(r)explicit} 
is absolutely convergent. Thus, assuming GRH, we have
\begin{equation} \label{zerosumbound} \sum_{\rho}\Gamma(\rho)N^{\rho} \ll N^{1/2}\sum_{\gamma } e^{-|\gamma|}\ll N^{1/2}\log 2 q. \end{equation}
Hence, for $\chi\neq \chi_0$,
\begin{equation} \label{Psi_r}
\Psi(r,\chi) \ll N^{1/2}\log 2q.
\end{equation}

\noindent\textit{\textbf{Proof} (Proof of \eqref{thm2a} in Theorem \ref{thm-smoothed})}
We note that \eqref{thm2a} is trivially true if $q> N^2$ since then
\[ F_q(N) \ll \sum_{\substack{n>N^2 \\ q|n}} n(\log n)^2 e^{-n/N} \ll N^2 e^{-N}\sum_{m} m^2e^{-m} \ll 1,\] and the error term in \eqref{thm2a} is larger than both main terms. Thus we assume henceforth that $2\le q \le N^2$, so that, for example, we can say $1\ll \log q \ll \log N$. Applying \eqref{Fqformula} and \eqref{Psi_r}, we have
\[F_q(N) = \frac{1}{\phi(q)}\Psi(r,\chi_0)^2 + O\left(\frac{1}{\phi(q)}\sum_{\chi \neq \chi_0}N\log^2\!q \right)+O\left((\log N\log q)^2\right) =\frac{1}{\phi(q)}\Psi(r,\chi_0)^2 + O\left(N \log^2\!q\right). \]

Therefore, by \eqref{PsiChi_0}, \eqref{Fq=1}, and the bound $\Psi(r) \ll N$, we have 
\[ \begin{split} F_q(N) &= \frac{1}{\phi(q)}(\Psi(r) + O(\log N \log q))^2 + O\left(N \log^2\!q\right) \\
&= \frac{1}{\phi(q)}\left(F(N) + O(N\log N\log q)\right) + O\left(N \log^2\!q)\right)\\
&= \frac{1}{\phi(q)}F(N) + O\left(N\log N\log q\right).
\end{split} \]
\qed

\noindent\textit{\textbf{Proof} (Proof of \eqref{thm2b} in Theorem \ref{thm-smoothed})} The idea here is that the smooth Fujii formula \eqref{thm2b} is just the square of the explicit formula for the power series generating function for primes. 
By \eqref{PsiNexplicit} in Lemma \ref{Psi(z)Explicit} we have
$\Psi(r) = N - \sum_{\rho}\Gamma(\rho) N^\rho +O(1)$.
Assuming RH, we have from \eqref{zerosumbound} that when $q=1$ that $\sum_{\rho}\Gamma(\rho) N^\rho \ll N^{1/2}$. Therefore by \eqref{Fq=1} on squaring we have
\[ F(N) = \Psi(r)^2 = \left(N - \sum_{\rho}\Gamma(\rho) N^{\rho} +O(1) \right)^2 
=N^2 -2\sum_{\rho}\Gamma(\rho) N^{\rho+1} + O(N).
\]

\qed

\section{ Transition from \texorpdfstring{$F_q(N)$}{Fq(N)} to \texorpdfstring{$G_q(N)$}{Gq(N)} }
\label{sec4}

In \cite{GS21} we made the transition from $F(z)=\Psi(z)^2$ to $G(N)$ to obtain \eqref{LZ}. A similar procedure works to make the transition from $F_q(z)$ to $G_q(N)$ to obtain \eqref{Gqformula}. Here we present only a sketch of the proof. The interested reader can fill in the details from \cite{GS21}.

The starting point is the formula, with $z=re(\alpha)$ and $r=e^{-1/N}$,
\begin{equation} \label{GqfromFq}
\int_0^1 F_q(z) I_N({\textstyle{\frac{1}{z}}})\, d\alpha = \sum_{\substack{n\\ q\mid n}}\psi_2(n)r^n 
\sum_{n'\le N} r^{-n'} \int_0^1 e(\alpha(n-n')) \, d\alpha = G_q(N),
\end{equation}
where
\begin{equation*} 
I_N(z) := \sum_{n\le N} z^n = z\left(\frac{1-z^N}{1-z}\right).
\end{equation*}
When $q=1$ we obtain
\begin{equation} \label{GfromPsi^2} G(N) = \int_0^1 \Psi(z)^2 I_N({\textstyle{\frac{1}{z}}})\, d\alpha . \end{equation}
We now take $2\le q\le N$. Applying Lemma \ref{Fqlem} to \eqref{GqfromFq} yields
\begin{align}
G_q(N) &= \frac{1}{\phi(q)}\int_0^1 \Psi(z,\chi_0)^2 I_N({\textstyle{\frac{1}{z}}})\, d\alpha
+ \frac{1}{\phi(q)}\sum_{\chi\neq \chi_0} \chi(-1) \int_0^1 \Psi(z,\chi)\Psi(z,\overline{\chi}) I_N({\textstyle{\frac{1}{z}}})\, d\alpha +
O\left(\log^5N\right) \notag\\
&= \frac{1}{\phi(q)}G(N) 
+ O\left(\max_{\chi\neq\chi_0} \int_0^{1/2} |\Psi(z,\chi)|^2 \min\{N, {\textstyle \frac{1}{\alpha}}\}\, d\alpha\right)
+ O\left(\frac{N}{\phi(q)}\log^3N\right) + O(\log^5N), \label{mainequation}
\end{align}
where we have used \eqref{PsiChi_0} and \eqref{GfromPsi^2}, while the first and third error terms are obtained by noting that
\[ \int_0^1 |I_N({\textstyle{\frac{1}{z}})}|\, d\alpha = \int_{-1/2}^{1/2} |I_N({\textstyle{\frac{1}{z}})}|\, d\alpha
\ll \int_{-1/2}^{1/2} \min\{N, {\textstyle \frac{1}{|\alpha|}}\}\, d\alpha
\ll \log N. \]

Following the proof of \cite[Theorem 2]{GS21}, we can bound the integral $\mathcal{I}$ in the first error term of \eqref{mainequation} as
\begin{equation} \label{calI} \mathcal{I} \ll \sum_{k\ll \log N}\frac{N}{2^k}\int_0^{\frac{2^k}{N}}|\Psi(z,\chi)|^2\, d\alpha, \end{equation}
where upon using Gallagher's lemma \cite{Montgomery71}, we have
\begin{equation} \label{key1} \int_0^{1/2h} |\Psi(z,\chi)|^2 \, d\alpha \ll \frac{1}{h^2}( I_1(N,h) + I_2(N,h)), \end{equation}
with
\begin{equation*} I_1(N,h) := \int_0^h \left|\sum_{n\le x}\chi(n)\Lambda(n)e^{-n/N}\right|^2\, dx
\end{equation*}
and
\begin{equation*}
I_2(N,h) :=\int_{0}^\infty \left| \sum_{x<n\le x+h}\chi(n)\Lambda(n)e^{-n/N}\right|^2\, dx. \end{equation*}
These mean values are closely related to the more familiar mean values, for $1\le h \le N$,
\begin{equation}\label{J1} J_1(h) := \int_0^h \left|\sum_{n\le x}\chi(n)\Lambda(n)\right|^2\, dx \end{equation}
and
\begin{equation}\label{J2}
J_2(N,h) :=\int_{0}^N \left| \sum_{x<n\le x+h}\chi(n)\Lambda(n)\right|^2\, dx. \end{equation}
By \cite[Lemma 2]{GV1996}, we have assuming GRH for $L(s,\chi)$ with $\chi$ primitive\footnote[3]{If $\chi$ in \eqref{J1} and \eqref{J2} is imprimitive one can still bound $J_1$ and $J_2$ by replacing $\chi$ with the primitive character $\chi^*$ which induces it in the sums in the integrands of \eqref{J1} and \eqref{J2} and add an insignificant error term of $O(\log N\log q)$ to the resulting sums.} that
\begin{equation*} J_1(h) \ll h^2\log^2\!q, \qquad \text{and} \qquad J_2(N,h) \ll hN\log^2\frac{3qN}{h}. \end{equation*}

Finally, using partial summation and Cauchy-Schwarz we find that the bounds for $J_1$ and $J_2$ also apply for $I_1$ and $I_2$, and therefore we obtain from \eqref{key1} on GRH that 
\begin{equation*} 
\int_0^{1/2h} |\Psi(z,\chi)|^2 \, d\alpha \ll \frac{N}{h}\log^2\!N.
\end{equation*}
Hence from \eqref{calI} we obtain $\mathcal{I} \ll N\log^3\!N$ and from \eqref{mainequation} we have, on GRH,
\begin{equation*} G_q(N) = \frac{1}{\phi(q)}G(N) +O(N\log^3\!N), \end{equation*}
which gives \eqref{Gqformula} with $C=3$.



\begin{thebibliography}{ABCD99}

\bibitem[BR18]{BhowRuzsa2018} G. Bhowmik, I. Z. Ruzsa, {\it Average Goldbach and the quasi-Riemann hypothesis}, Anal. Math. {\bf 44} (2018), no. 1, 51--56. 

\bibitem[BS10]{BhowPuchta2010} G. Bhowmik, J. -C. Schlage-Puchta, {\it Mean representation number of integers as the sum of primes}, Nagoya Math. J. {\bf 200} (2010), 27--33.

\bibitem[BHMS19]{Bhowmik-H-M-SGoldbach2019} Gautami Bhowmik, Karin Halupczok, Kohji Matsumoto, Yuta Suzuki, {\it Goldbach representations in arithmetic progressions and zeros of Dirichlet L-functions}, Mathematika {\bf 65} (2019), no. 1, 57--97.
 





\bibitem[FGIS22]{FGIS21} J. B. Friedlander, D. A. Goldston, H. Iwaniec, and A. I. Suriajaya, {\it Exceptional zeros and the Goldbach Problem}, J. Number Theory {\bf 233} (2022), 78--86.

\bibitem[Fuj91a]{Fujii1} A. Fujii, {\it An additive problem of prime numbers}, Acta Arith. {\bf 58} (1991), 173--179.

\bibitem[Fuj91b]{Fujii2} A. Fujii, {\it An additive problem of prime numbers. II}, Proc. Japan Acad. ser. A Math. Sci. {\bf 67} (1991), 248--252.

\bibitem[Fuj91b]{Fujii3} A. Fujii, {\it An additive problem of prime numbers. III}, Proc. Japan Acad. ser. A Math. Sci. {\bf 67} (1991), 278--283.

\bibitem[GS23]{GS21} D. A. Goldston and A. I. Suriajaya, {\it On an average Goldbach representation formula of Fujii}, to appear in Nagoya Math. J., preprint in arXiv:2110.14250 [math.NT].

\bibitem[GV96]{GV1996} D. A. Goldston and R. C. Vaughan, {\it On the Montgomery-Hooley asymptotic formula}, Sieve Methods, Exponential Sums, and their Application in Number Theory, Greaves, Harman, Huxley Eds., Cambridge University Press, 1996, 117--142.

\bibitem[GY17]{Gold-Yang} D. A. Goldston and L. Yang, {\sl The Average Number of Goldbach Representations}, in: Prime numbers and representation theory, Edited by Ye Tian \& Yangbo Ye, Science Press, Beijing, 2017, 1--12.

\bibitem[Gran07]{Gran07} A. Granville, {\it Refinements of Goldbach's conjecture, and the generalized Riemann hypothesis}, Funct. Approx. Comment. Math. {\bf 37} (2007), 159--173.

\bibitem[Gran08]{Gran08} A. Granville, {\it Corrigendum to \lq\lq Refinements of Goldbach's conjecture, and the generalized Riemann hypothesis\rq\rq}, Funct. Approx. Comment. Math. {\bf 38} (2008), 235--237.

\bibitem[HL18]{HardyLittlewood1918} G. H. Hardy and J. E. Littlewood, {\it Contributions to the theory of the Riemann zeta-function and the theory of the distribution of primes}, Acta Math. {\bf 41} (1918), 119--196. Reprinted as pp. 20--97 in {\it Collected Papers of G. H. Hardy}, Vol. II, Clarendon Press, Oxford University Press, Oxford, 1967.

\bibitem[HL19]{HardyLittlewood1919} G. H. Hardy and J. E. Littlewood, {\it Note on Messrs. Shah and Wilson's paper entitled: On an empirical formula connected with Goldbach's Theorem}, Proceedings of the Cambridge Philosophical Society, vol. {\bf 19} (1919), 245--254. Reprinted as pp. 535--544 in {\it Collected Papers of G. H. Hardy}, Vol. I, Clarendon Press, Oxford University Press, Oxford, 1966.

\bibitem[HL22]{HardyLittlewood1922} G. H. Hardy and J. E. Littlewood, {\it Some problems of `Partitio numerorum'; III: On the expression of a number as a sum of primes}, Acta Math. {\bf 44} (1922), no. 1, 1--70. Reprinted as pp. 561--630 in {\it Collected Papers of G. H. Hardy}, Vol. I, Clarendon Press, Oxford University Press, Oxford, 1966.


\bibitem[Ing32]{Ingham1932} A. E. Ingham, {\it The Distribution of Prime Numbers}, Cambridge Tracts in Mathematics and Mathematical Physics 30, Cambridge Univ. Press, Cambridge, 1932.


\bibitem[Lan00]{Landau1900} E. Landau, {\it Ueber die zahlentheoretische Funktion $\phi(n)$ und ihre Beziehung zum Goldbachschen Satz}, G{\"o}ttinger Nachrichten, 1900, 177--186.

\bibitem[LZ12a]{Lang-Zac1} A. Languasco and A. Zaccagnini, {\it The number of Goldbach representations of an integer}, Proc. Amer. Math. Soc. {\bf 140} (2012), 795--804.

\bibitem[LZ12b]{Lang-Zac2} A. Languasco and A. Zaccagnini, {\it Sums of many primes}, Journal of Number Theory {\bf 132} (2012), 1265--1283.

\bibitem[Mon71]{Montgomery71} Hugh L. Montgomery, {\it Topics in Multiplicative Number Theory}, Lecture Notes in Mathematics, Vol. 227, Springer-Verlag, Berlin-New York, 1971.


\bibitem[MV07]{MontgomeryVaughan2007} H. L. Montgomery and R. C. Vaughan, {\it Multiplicative Number Theory}, Cambridge Studies in Advanced Mathematics {\bf 97}, Cambridge University Press, Cambridge, 2007.

\bibitem[Suz17]{Suzuki2017} Y. Suzuki, {\it A mean value of the representation function for the sum of two primes in arithmetic progressions}, Int. J. Number Theory {\bf 13} (4) (2017), 977--990.

\end{thebibliography}
\end{document}